\documentclass[12pt,leqno,fleqn,epsfig]{article}
\usepackage{amssymb, epsfig, amsmath, amsthm}
\usepackage{mathrsfs}       
\usepackage{color}  
     
\numberwithin{equation}{section}

\newtheorem{thm}{Theorem}[section]

\newtheorem{cor}{Corollary}[section]

\renewcommand{\t}{\theta}
\renewcommand{\t}{\theta}

\newcommand{\g}{\gamma}

\renewcommand{\a}{\alpha}
\renewcommand{\b}{\beta}

\renewcommand{\Omega}{\Theta}

\newcommand{\bb}{\begin{equation}}
\newcommand{\ee}{\end{equation}}
\newcommand{\bq}{\begin{eqnarray}}
\newcommand{\eq}{\end{eqnarray}}
\newcommand{\bqn}{\begin{eqnarray*}}
\newcommand{\eqn}{\end{eqnarray*}}

\newcommand{\nrm}[1]{\Vert#1\Vert}

\newcommand{\nnrm}[1]{{\vert\kern-0.25ex\vert\kern-0.25ex\vert #1 
		\vert\kern-0.25ex\vert\kern-0.25ex\vert}}

\newcommand{\rd}{\partial}

\newcommand{\gmm}{\gamma}

\newcommand{\omg}{\omega}

\newcommand{\bbR}{\mathbb R}

\begin{document}
\title{Preservation of  log-H\"{o}lder coefficients of the vorticity in the transport equation}

\author{Dongho Chae\thanks{Department of Mathematics, Chung-Ang University.  E-mail: dchae@cau.ac.kr} \and In-Jee Jeong\thanks{Department of Mathematical Sciences and RIM, Seoul National University.  E-mail: injee\_j@snu.ac.kr}}
\date{\today}
\maketitle

\begin{abstract}
We show that the log-H\"{o}lder coefficients of a solution to the transport equation is  preserved in time.
\end{abstract}
\section{Introduction}
\setcounter{equation}{0}
We are concerned on  the transport equations in $ \Bbb R^n$:
\bb\label{tr}
\t_t+u\cdot \nabla\theta= 0 \quad\text{on}\quad \Bbb R^n \times (0, T),
\ee
 where $u=u(x,t)=(u_1(x,t), \cdots, u_n (x,t))$ for $(x,t) \in \Bbb R^n \times [0, +\infty)$, and $\theta =\theta (x,t)$ is a scalar
 with $\theta(x,0)=\theta_0(x)$.\\

 We say $f$ is H\"{o}lder continuous with exponent $\g\in (0, 1)$ at $x_0\in \Bbb R^n$ 
 if there exists $ B_r(x_0) \subset \Bbb R^n$ such that 
 $$
  \sup_{ x\in B_r(x_0)} \frac{ |f (x) -f(x_0) |}{ |x-x_0|^\g}<+\infty.
 $$
 In this case we define the H\"{o}lder coefficient  $[f]_{\g;x_0}$  by 
  \bb\label{co1a}
[f]_{\g;x_0}  =   \lim_{r\to 0} \sup_{ x\in B_r (x_0)} \frac{ |f (x) -f(x_0) |}{ |x-x_0|^\g},
\ee 
 and denote 
$$ C^\g (x_0)=\{ f\in L^\infty (\Bbb R^n )\, |\, [f]_{ \g; x_0 } <+\infty \}.
$$
Note that we are following  the notation of \cite[pp.52]{gil}. 
One can generalize the above notion to a general  case of functions defining the modulus of continuity $\delta=\delta (s)$,  and denote
  \bb\label{co1}
 [f]_{\delta(s);x_0}  :=   \lim_{r\to 0} \sup_{ x\in B_r (x_0)}\frac{ |f (x) -f(x_0) |}{\delta ( |x-x_0| )},
\ee  
for $x_0\in \Bbb R^n$, and write  
$$ C^{\delta(s) } (x_0 )=  \{ f\in L^\infty (\Bbb R^n) \, |\, [f]_{\delta(s);x_0}  <+\infty \}.
$$
Given a vector field $u=u(x,t)\in L^1 (0, T; Lip (\Bbb R^n))$, we introduce the particle trajectory mapping  $x_0 \mapsto \varphi(x_0 , t)$ from $\Bbb R^n$ into $\Bbb R^n$   generated by $u=u(x,t)$, which  is defined as the solution of  the ordinary differential equation,
\bb\label{particle}
\left\{ \aligned &\frac{\partial   \varphi(x_0, t)}{\partial t}   =u( \varphi (x_0, t)  , t) \quad \text{on}\quad (0,T),\\
                    &\varphi(x_0 , 0)  = x_0 \in \Bbb R^n.\endaligned
                    \right.
\ee

Our main result states that for the transport equation, the log-H\"older coefficient is constant in time. 
 \begin{thm}
Let    $u\in L^1 (0, T; Lip (\Bbb R^n))$ and $\theta \in L^\infty(0,T; L^\infty(\Bbb R^n))$ satisfies \eqref{tr}  in $\Bbb R^n\times (0, T)$.
 Suppose there exists $ \g \in (0,\infty )$ such that $\t_0 \in C^{(\log \frac{1}{s} )^{-\g}} (x_0)$ for some  $x_0 \in \Bbb R^n$.
Then, 
\bb\label{th2aa}
[\t(\cdot,  t) ]_{ {\log(\frac{1}{s})}^{-\g}; \varphi (x_0, t) } =
 [\t_0 ]_{ {\log(\frac{1}{s})}^{-\g}; x_0} \quad  \forall t\in [0, T),
\ee
where $\varphi(x_0, t)$ denotes the particle trajectory defined in \eqref{particle}. In particular, we have 
\bb\label{th2b}
\sup_{x\in \Bbb R^n} [\theta (t) ]_{ {\log(\frac{1}{s})}^{-\g}, x } =
\sup_{x\in \Bbb R^n} [\theta_0 ]_{ {\log(\frac{1}{s})}^{-\g}, x } .
\ee
\end{thm}

As an application we consider  the 2D Euler equations on $\Bbb R^2\times [0, +\infty)$. 
$$
(E) \left\{\aligned  &\omega_t+u\cdot \nabla \omega= 0,\label{e1}\\
& \nabla \cdot u=0,\quad \partial_1 u_2-\partial_2 u_1=\omega,\\
& \omega(x,0)=\omega_0(x).
\endaligned \right.
$$
For bounded initial data $\omg_0 \in L^1 \cap L^\infty (\bbR^2)$, well known Yudovich theory gives a unique solution to $(E)$ in $\omg(t,x) \in L^\infty_t (L^1 \cap L^\infty)_x.$ If $\omg_0$ has a little bit of additional regularity, then it is guaranteed that $u$ is Lipschitz continuous in space (\cite{MB}).

  \begin{cor}
Let    $u\in L^1 (0, T; Lip (\Bbb R^2))$  be a solution to  (E), 
Suppose $\omg_0 \in C^{(\log \frac{1}{s} )^{-\g}} (x_0)$ for some  $ \g \in (0,\infty )$ and $x_0 \in \Bbb R^2$.
Then, 
\bb\label{th2}
[\omega (\cdot,  t) ]_{ {\log(\frac{1}{s})}^{-\g}; \varphi (x_0, t) } =
 [\omega_0 ]_{ {\log(\frac{1}{s})}^{-\g}; x_0} \quad  \forall t\in [0, T).
\ee
\end{cor}

{\bf Remark 1.1 } In the proof, we observe the following inequalities.
\begin{align}\label{th1}
&[\omega_0 ]_{ \beta; x_0}\exp \left(-\beta \int_0 ^t  |\nabla u (\varphi (x_0,s),s) | ds\right)\cr
&\qquad \le [\omega (\cdot,  t) ]_{ \beta; \varphi (x_0, t) } \le [\omega_0 ]_{ \beta; x_0}\exp \left(\beta \int_0 ^t   |\nabla u (\varphi (x_0,s),s) | ds\right) \quad \forall \beta \in (0, 1]
\end{align}
for all $t\in [0, T)$. We could understand that \eqref{th2} is a limiting relation of \eqref{th1} as $\beta \to 0$.\\
\ \\
{\bf Remark 1.2 }   It would be interesting to compare \eqref{th2}  with   the result of the double exponential growth in time of the vorticity gradient by Kieselev and \v{S}ver\'{a}k\cite{kie} for the 2D Euler equations. As far as the authors know we do not know if the growth of $\sup_{x\in \Bbb R^2} [\omg (t)]_{\g, x}$ (or $\| \omg (t)\|_{C^{0, \g }}$)  is slower than the double exponential one. The above result says that for the log-H\"{o}lder  coefficient it is constant in time.

\ \\
{\bf Remark 1.3 } In the 2D Euler equations, initial vorticity with asymptotics \begin{equation*}
	\begin{split}
		|\omg_0(x)| \simeq \log(\frac{1}{|x|})^{-\gmm}, \qquad 0< |x| \ll 1
	\end{split}
\end{equation*} for some $\gmm>0$ together with odd-odd symmetry with respect to both axes and satisfying $\omg_0 \ge 0$ on $(\bbR_+)^2$ was considered in \cite{EJ}.\footnote{This type of data could be used to prove \textit{ill-posedness} of the 2D Euler equations in the critical Sobolev space $\omg\in H^1$.} Our main result shows that whenever $\gmm\le1$, \textit{if} the corresponding velocity field belongs to $L^1(0,T;Lip)$ for any $T>0$, it is forced that $|\omg(x,t)| \simeq \log(\frac{1}{|x|})^{-\gmm}$. Here, the origin is fixed by the flow due to the odd symmetry. But then, using this asymptotics with the explicit formula \begin{equation*}
\begin{split}
	(\rd_{x_1} u_1)(0,t) = c_0 \int_{(\bbR_+)^2} \frac{y_1y_2}{|y|^4} \omg(y,t) \, dy \gtrsim \int_{r(t)}^{\infty} \frac{1}{r (\log\frac{1}{r})^\gmm} \, dr , 
\end{split}
\end{equation*} one can prove that $|\nabla u(\cdot,t)| = +\infty$ for $t\in (0,T)$, which is a contradiction to the $L^1(0,T;Lip)$--assumption. Therefore, as a corollary of our main result, we deduce that for such an initial vorticity, $\int_0^T  \nrm{ \nabla u(\cdot, t)}_{L^\infty} \, dt = +\infty$ for any $T>0$. In turn, this forces the \textit{Lagrangian deformation} to be divergent; $|\rd_{x_1} \Phi_1(0,t)|=+\infty$ for $t>0$. This is very interesting, since a formal analysis given in \cite{EJSVP1} suggests that $\nrm{ \nabla u(\cdot, t)}_{L^\infty} \simeq ct^{-1}$ for some absolute constant $c>0$, which is barely non-integrable in time. To clarify, in general (without the odd-odd assumption) it is possible that $|\omg(x,t)| \simeq \log(\frac{1}{|x|})^{-\gmm}$ and $ u(x,t) \in L^1_t Lip$ hold at the same time for any $\gmm>0$. This happens for instance the vorticity satisfies a rotational symmetry; see \cite{EJS}.

\ \\ 

For the three-dimensional incompressible Euler equations, we have the following \begin{thm}
	Let $(u,\omega)$ be a solution to the 3D Euler equations on $\mathbb{R}^3 \times [0,T)$, namely \begin{equation}\label{eq:3D-Euler}
		\left\{
		\begin{aligned}
			&\partial_t \omg + u \cdot \nabla \omg = \omg \cdot \nabla u, \\
			&\nabla \cdot u = 0, \quad \nabla \times u = \omg. 
		\end{aligned}
		\right.
	\end{equation}  Let $\g >0$, and suppose there exists $\g_1 >\g $ such that
	$$
	\int_0 ^t   [u(\tau)]_{1; \varphi(x_0, \tau) }  [ \nabla u (\tau ) ]_ { {\log(\frac{1}{s})}^{-\g_1};\varphi(x_0, \tau ) }  d\tau <+\infty
	$$ holds for some $0<t <T$. 
	Then,  we have 
	\bb\label{th2aaa}
	[\omg(\cdot,  t) ]_{ {\log(\frac{1}{s})}^{-\g};\varphi(x_0, t) } =
	[\omg_0 ]_{ {\log(\frac{1}{s})}^{-\g}; x_0} \quad \forall x_0 \in \Bbb R^n.
	\ee
	In particular,
	\bb
	\sup_{x\in \Bbb R^n} [\omg(\cdot,  t) ]_{ {\log(\frac{1}{s})}^{-\g}; x } =
	\sup_{x\in \Bbb R^n}[\omg_0 ]_{ {\log(\frac{1}{s})}^{-\g}; x}.
	\ee
\end{thm}  
  
\section{ The Proof of the Main Theorem}
\setcounter{equation}{0}  
    \noindent{\bf Proof  of Theorem 1.1 } In the proof, we assume that $u(x,t)$ is $C^1$--smooth in $x$, so that the flow map $\varphi(\cdot,t)$ defines a $C^1$--diffeomorphism on $\bbR^2$ for $0\le t<T$. The case of $u(x,t) \in L^1(0,T;Lip(\bbR^n))$ then follows from a simple approximation argument. 
    
    From $\varphi(\a, t)= \a +\int_0 ^t u(\varphi(\a, t),t) ds$ and the triangle inequality, 
    \begin{align}
    	& |\a -\b | -\left|\int_0 ^t  \left( u(\varphi(\a, t),t) - u(\varphi(\b, t),t) \right)  ds\right|\cr
    	& \qquad\le |\varphi(\a,t) -\varphi(\b,t) | \le |\a -\b | +\left|\int_0 ^t  \left( u(\varphi(\a, t),t) - u(\varphi(\b, t),t) \right)  ds\right|\end{align}
    from which, diving  it by $|\a-\b|\not=0$,  we have
    \begin{align}
    	&1- \int_0 ^t \frac{ |u (\varphi(\a, s),s )-u (\varphi(\b, s), s) |}{ |\varphi(\a, s)-\varphi(\b, s)| } \frac{ |\varphi(\a,s)-\varphi(\b, s) |}{|\a -\b|} ds\cr
    	&\qquad\le \frac{ |\varphi(\a,t)-\varphi(\b, t) |}{|\a -\b|} \cr
    	&\qquad\le 1+  \int_0 ^t \frac{ |u (\varphi(\a, s),s )-u (\varphi(\b, s), s) |}{ |\varphi(\a, s)-\varphi(\b, s)| } \frac{ |\varphi(\a,s)-\varphi(\b, s) |}{|\a -\b|} ds,
    \end{align}
    and 
    \begin{align}
    	&1- \int_0 ^t \|\nabla u(s)\|_{L^\infty}  \frac{ |\varphi(\a,s)-\varphi(\b, s) |}{|\a -\b|} ds\cr
    	&\qquad\qquad\le \frac{ |\varphi(\a,t)-\varphi(\b, t) |}{|\a -\b|} \le 1+  \int_0 ^t \|\nabla u(s)\|_{L^\infty} \frac{ |\varphi(\a,s)-\varphi(\b, s) |}{|\a -\b|} ds.
    \end{align}
    By Gronwall's inequality we obtain
    \bb\label{key}
    \frac{1}{\mu(t)} \le  \frac{ |\varphi(\a,t)-\varphi(\b, t) |}{|\a -\b|} 
    \le \mu(t), \qquad \mu(t)=e ^ { \int_0 ^t  \|\nabla u(s)\|_{L^\infty}  ds}
    \ee
    for all $\a, \b \in \Bbb R^n$.
    This shows that
    $$
    \b\in B_{ \frac{r}{\mu(t)}} (\a ) \Rightarrow \varphi(\b,t) \in B_r ( \varphi(\a,t)) \Rightarrow \b\in B_{ \mu(t) r} (\a)\quad \forall r>0,
    $$
    which implies that 
    \bb
    \varphi( B_r (\a), t) \subset B_{\mu(t) r} ( \varphi(\a,t ))\quad \text{and}\quad B_r (\varphi(\a, t)) \subset  \varphi (B_{\mu(t) r} (\a ), t).
    \ee
    From this we infer
    \bb\label{key1}
    \sup_{\b \in B_r (\a)} f( \varphi(\b, t) )=  \sup_{\varphi(\b, t)  \in \varphi( B_r (\a), t) } f( \varphi(\b, t) )   \le \sup_{ y \in B_{ \mu(t) r} ( \varphi(\a, t) ) } f( y ),
    \ee
    and
    \bb\label{key2}
    \sup_{y \in B_r (\varphi(\a,t))} f( y) \le  \sup_{y   \in \varphi( B_{\mu (t)r}  (\a), t) } f( y ) = \sup_{ \b \in B_{ \mu(t) r} (  \a) } f( \varphi(\b, t) ).
    \ee
    Combining \eqref{key1} with \eqref{key2}, we have
    \begin{align*}
    	\lim_{r\to 0} \sup_{\b \in B_r (\a)} f( \varphi(\b, t) )&\le \lim_{r\to 0} \sup_{ y \in B_{ \mu(t) r} ( \varphi(\a, t) ) } f( y ) = \lim_{r\to 0} \sup_{y \in B_r (\varphi(\a,t))} f( y)\cr
    	&\le \lim_{r\to 0} \sup_{ \b \in B_{ \mu(t) r} (  \a) } f( \varphi(\b, t) ) =  \lim_{r\to 0} \sup_{\b \in B_r (\a)} f( \varphi(\b, t) ).
    \end{align*}
    Hence,
    \bb\label{key3}
    \lim_{r\to 0} \sup_{\b \in B_r (\a)} f( \varphi(\b, t) )=  \lim_{r\to 0} \sup_{y \in B_r (\varphi(\a,t))} f( y).
    \ee
    From \eqref{key} one has
    \begin{align} 
    	\left( 1+ \frac{ \log \mu(t) }{ \log \frac{1}{|\a-\b|} }  \right) ^{-\g}   \le  \frac{ \left(\log  \frac{1}{ |\varphi(\a,t)-\varphi(\b, t) |}\right) ^{-\g}}{\left( \log \frac{1}{|\a-\b|} \right) ^{-\g}} \le \left( 1-\frac{\log \mu(t) }{\log \frac{1}{|\a-\b|} }\right) ^{-\g},
    \end{align} 
    and thus
    \bb\label{key4}
    \lim_{r\to 0} \sup_{\b\in B_r(\a)} \frac{ \left(\log  \frac{1}{ |\varphi(\a,t)-\varphi(\b, t) |}\right) ^{-\g}}{\left( \log \frac{1}{|\a-\b|} \right) ^{-\g}} =1.
    \ee
    From  $\theta (\varphi(\a, t), t)=\theta_0 (\a )$ we deduce
    \begin{align} \label{tr}
    	[\theta_0]_{(\log\frac{1}{s})^{-\g}, \a} &=\lim_{r\to 0} \sup_{\b\in B_r(\a)} \frac{ |\theta_0 (\a )- \theta_0 (\b)|}{ \left(\log \frac{1}{|\a-\b|} \right)^{-\g} } \cr
    	&=  \lim_{r\to 0} \sup_{\b\in B_r(\a)}\frac{ |\theta (\varphi(\a, t), t) -\theta (\varphi(\b, t),t)|}{ \left(\log\frac{1}{|\varphi(\a,t)-\varphi(\b, t) |} \right)^{-\g} }
    	\frac{ \left(\log\frac{1}{|\varphi(\a,t)-\varphi(\b, t) |} \right)^{-\g} }{ \left(\log \frac{1}{|\a-\b|} \right)^{-\g} }\cr
    	&= \lim_{r\to 0} \sup_{\b  \in B_r (\a)}\frac{ |\theta (\varphi(\a, t), t) -\theta (\varphi(\b,t) ,t)|}{ \left(\log\frac{1}{|\varphi(\a,t)-\varphi(\b,t) |} \right)^{-\g} }\cr
    	&
    	=\lim_{r\to 0} \sup_{y\in B_{ r} (\varphi(\a,t))}\frac{ |\theta (\varphi(\a, t), t) -\theta (y,t)|}{ \left(\log\frac{1}{|\varphi(\a,t)-y |} \right)^{-\g} }\cr
    	&=     [\theta(t)]_{(\log\frac{1}{s})^{-\g}, \varphi(\a, t)},
    \end{align}
    where we used \eqref{key4} in the third equality and \eqref{key3} in the fourth equality  respectively. $\square$\\
    \ \\

    {\bf Remark 2.1 }  The result is sharp in the case of the 2D Euler equations in the sense that for any $\beta>0$, there exists initial data $\omega_0 \in C^{\beta}(x_{0}) $ with $[\omega(\cdot,t)]_{\beta; \varphi(x_{0},t)}$ growing at least exponentially in time.  To see this, we take $\omega_{0}(x) \in C^{\beta}_c(\mathbb{R}^2)$ which is odd with respect to both axes, non-negative on the first quadrant and satisfies $$ \omega_{0}(x) = \frac{2x_{1}x_{2}}{4x_{1}^2 + x_{2}^{2}} |x|^{\beta}$$ for $|x|\le 1$. There is a unique global in time solution (with Lipschitz continuous velocity) for any $\beta>0$. From the symmetry, $\varphi(0,t)=0$ for all $t$. For fixed $|x|\le 1$, the maximum of $\omega_{0}(x)$ is achieved for $2x_1 = x_2$. Using the version of Kiselev-Sverak ``Key Lemma'' appeared in Zlatos \cite{Z}, one can prove at least for a short time interval that \begin{equation*}
    	\begin{split}
    		\varphi_1(x,t) \simeq r \exp\left( (\frac{C_{0}}{\beta} + O(1))t  \right), \qquad \varphi_2(x,t) \simeq 2r \exp\left( -(\frac{C_{0}}{\beta} + O(1))t  \right),
    	\end{split}
    \end{equation*} for all $x = (r,2r)$ with $0<r<1/10$. In particular, $|\varphi(x,t)|$ is decreasing in time. Then, using $$\omega(\varphi(x,t),t) = \omega_{0}(x)= \frac{2x_{1}x_{2}}{4x_{1}^2 + x_{2}^{2}} |x|^{\beta} = \frac{1}{2}(\sqrt{5}r)^{\beta}$$ for $x = (r,2r)$ and \begin{equation*}
    	\begin{split}
    		[\omega(\cdot,t)]_{\beta; 0} \ge \lim_{r\to 0} \frac{ (\sqrt{5}r)^{\beta}  }{ 2|\varphi(x,t)|^{\beta} } , \qquad x = (r,2r), 
    	\end{split}
    \end{equation*} we can show local in time growth of $[\omega(\cdot,t)]_{\beta; 0}$. This argument could work to prove growth of moduli of continuity which lie ``between'' H\"older and log-H\"older.
    
    \ \\
    
     \noindent{\bf Proof  of Theorem 1.3 }
    From the vorticity form of the Euler equations
    $$
    \omg_t +u\cdot \nabla \omg =\omg\cdot \nabla u, \quad \nabla \times u=\omg, \quad \nabla \cdot u=0,
    $$
    we have
    $$
    \omg (\varphi(\a, t), t)= \omg_0 (\a) +\int_0 ^t (\omg \cdot \nabla) u ( \varphi(\a, \tau), \tau ) d\tau.
    $$
    By the triangular inequality
    \begin{align}\label{pr201}
    	&\frac{|\omg_0 (\a)-\omg_0 (\b)|}{ \left(\log\frac{1}{|\a-\b|} \right)^{-\g}} -\int_0 ^t\frac{|\omg\cdot \nabla u (\varphi(\a,\tau),\tau)-\omg\cdot \nabla u(\varphi(\b, \tau), \tau) |}{ \left(\log\frac{1}{|\a-\b|} \right)^{-\g}} d\tau\cr
    	&\qquad\le\frac{|\omg(\varphi(\a, t),t)-\omg(\varphi(\b, t),t)|}{ \left(\log\frac{1}{|\a-\b|} \right)^{-\g}} \cr
    	&\qquad \le  \frac{|\omg_0 (\a)-\omg_0 (\b)|}{ \left(\log\frac{1}{|\a-\b|} \right)^{-\g}}
    	+\int_0 ^t\frac{| \omg\cdot \nabla u(\varphi(\a,\tau),\tau)-\omg\cdot \nabla u(\varphi(\b, \tau), \tau) |} { \left(\log\frac{1}{|\a-\b|} \right)^{-\g}} d\tau
    \end{align} 
    We note that
    \begin{align}
    	&| \omg\cdot \nabla u(\varphi(\a,\tau),\tau)-\omg\cdot \nabla u(\varphi(\b, \tau), \tau) |
    	\\
    	&\le |\omg (\varphi(\a, \tau), \tau)| |\nabla u (\varphi(\a, \tau), \tau)- \nabla u (\varphi(\b, \tau), \tau)|\cr
    	&\qquad+ |\nabla u (\varphi(\b, \tau), \tau)| |\omg (\varphi(\a, \tau), \tau) - \omg (\varphi(\b, \tau), \tau) |\cr
    	&\le 2|\nabla u(\varphi(\a, \tau), \tau)| |\nabla u (\varphi(\a, \tau), \tau)- \nabla u (\varphi(\b, \tau), \tau)|\cr
    	&\qquad+2 |\nabla u (\varphi(\b, \tau), \tau)| |\nabla u (\varphi(\a, \tau), \tau) - \nabla u(\varphi(\b, \tau), \tau) | \cr
    	&\le 2  (|\nabla u(\varphi(\a, \tau), \tau)|   + |\nabla u (\varphi(\b, \tau), \tau)| ) ( |\nabla u (\varphi(\a, \tau), \tau) - \nabla u (\varphi(\b, \tau), \tau) |). 
    \end{align}
    Hence, we have
    \begin{align}\label{vorr}
    	\hspace{-2.in}& \lim_{r\to 0} \sup_{\b \in B_r ( \a)}\int_0 ^t\frac{| \omg\cdot \nabla u(\varphi(\a,\tau),\tau)-\omg\cdot \nabla u(\varphi(\b, \tau), \tau) |} { \left(\log\frac{1}{|\a-\b|} \right)^{-\g}} d\tau \cr
    	&\le 4\lim_{r\to 0} \sup_{\b \in B_r ( \a)}\int_0 ^t   |\nabla u(\varphi(\a, \tau), \tau)|   \frac{( |\nabla u (\varphi(\a, \tau), \tau) - \nabla u (\varphi(\b, \tau), \tau) |)} { \left(\log\frac{1}{|\a-\b|} \right)^{-\g}} d\tau \cr
    	& \le 4\lim_{r\to 0} \sup_{\b \in B_r ( \a)}\int_0 ^t   |\nabla u(\varphi(\a, \tau), \tau)| \times \cr
    	&\qquad\qquad \times  \frac{( |\nabla u (\varphi(\a, \tau), \tau) - \nabla u (\varphi(\b, \tau), \tau) |)}{\left(\log\frac{1}{ |\varphi( \a, s)-\varphi(\b, \tau) |} \right)^{-\g_1}} \frac{\left(\log\frac{1}{ |\varphi( \a, \tau)-\varphi(\b, \tau) |} \right)^{-\g_1}} { \left(\log\frac{1}{|\a-\b|} \right)^{-\g}} d\tau\cr
    	&  \le4\int_0 ^t [u(t)]_{1; \varphi(\a, \tau )} [\nabla u]_{ (\log\frac{1}{s})^{-\g_1}; \varphi(\a, \tau)} \lim_{r\to 0} \sup_{\b \in B_r ( \a)}  \frac{\left(\log\frac{1}{ |\varphi( \a, \tau)-\varphi(\b, \tau) |} \right)^{-\g_1}} { \left(\log\frac{1}{|\a-\b|} \right)^{-\g}} d\tau\cr
    	&=0.
    \end{align}
    Taking $\lim_{r\to 0} \sup_{\b \in B_r ( \a)}$ on \eqref{pr201}, and substituting \eqref{vorr}, 
    $$\lim_{r\to 0} \sup_{\b \in B_r ( \a)}\frac{|\omg(\varphi(\a, t),t)-\omg(\varphi(\b, t),t)|}{ \left(\log\frac{1}{|\a-\b|} \right)^{-\g}} 
    = [\omg (t) ]_{ (\log\frac{1}{s})^{-\g_1}; \varphi(\a, \tau)},
    $$
    we obtain \eqref{th2aaa}. $\square$\\

\subsection*{Acknowledgments}{D.~Chae was supported partially  by NRF grant 2021R1A2C1003234. I.-J.~Jeong was supported by  the New Faculty Startup Fund from Seoul National University, the Science Fellowship of POSCO TJ Park Foundation, and the National Research Foundation of Korea grant (No. 2019R1F1A1058486).}

\end{document}